\documentclass[11pt,reqno]{article}
\usepackage{amsmath,mathrsfs,graphicx,amssymb,amsfonts,color}
\usepackage{cite}

\font\bg=cmbx10 scaled\magstep1

\font\small=cmr8
\newtheorem{newlemma}{{\bf Lemma}}

\newenvironment{lema}{\begin{newlemma}{\hspace{-0.5
em}{\bf.}}}{\end{newlemma}}

\newtheorem{newteorem}{{\bf Theorem}}

\newenvironment{teorem}{\begin{newteorem}{\hspace{-0.5
em}{\bf.}}}{\end{newteorem}}

\newtheorem{newkorolari}{{\bf Corollary}}

\newenvironment{korolari}{\begin{newkorolari}{\hspace{-0.5
em}{\bf.}}}{\end{newkorolari}}

\newtheorem{newdefine}{{\bf Definition}}

\newtheorem{newquestion}{{\bf Question}}

\newtheorem{newkonjek}{{\bf Conjecture}}

\newenvironment{konjek}{\begin{newkonjek}{\hspace{-0.5
em}{\bf.}}}{\end{newkonjek}}

\newtheorem{newexample}{{\bf Example}}

\begin{document}
\tolerance=10000
\baselineskip18truept
\newbox\thebox
\global\setbox\thebox=\vbox to 0.2truecm{\hsize 0.15truecm\noindent\hfill}
\def\boxit#1{\vbox{\hrule\hbox{\vrule\kern0pt
\vbox{\kern0pt#1\kern0pt}\kern0pt\vrule}\hrule}}
\def\qed{\lower0.1cm\hbox{\noindent \boxit{\copy\thebox}}\bigskip}
\def\nt{\noindent}

\centerline{\Large \bf More on energy and Randi\'{c} energy of specific  graphs }
\vspace{.3cm}

%\centerline {\Large \bf  are unimodal}
\bigskip

\baselineskip12truept
\centerline{Saeid  Alikhani$^{}${}\footnote{\baselineskip12truept\it\small
Corresponding author. E-mail: alikhani@yazd.ac.ir}and  Nima Ghanbari$^{}${}\footnote{\baselineskip12truept\it\small
E-mail: n.ghanbari@stu.yazd.ac.ir} }
\baselineskip20truept
\centerline{\it Department of Mathematics, Yazd University}
\vskip-8truept
\centerline{\it  89195-741, Yazd, Iran}

\vskip-0.2truecm
\noindent\rule{16cm}{0.1mm}
\noindent{\bg{Abstract}}

\baselineskip14truept

{\nt Let $G$ be a simple graph of order $n$. The energy $E(G)$ of the graph $G$ is the sum of the absolute values of the
eigenvalues of $G$. The Randi\'{c} matrix of $G$, denoted by $R(G)$, is defined as the $n\times n$ matrix whose $(i,j)$-entry is $(d_id_j)^{\frac{-1}{2}}$ if $v_i$ and
$v_j$ are adjacent and $0$ for another cases. The Randi\'{c} energy $RE$ of $G$ is the sum of
absolute values of the eigenvalues of $R(G)$. In this paper we compute the energy and Randi\'{c} energy for  certain graphs. Also we propose a conjecture on Randi\'c energy. }


\noindent\rule{16cm}{0.1mm}

\baselineskip20truept

\section{Introduction}

\nt In this paper we are concerned with simple finite graphs, without directed, multiple, or weighted edges, and without self-loops. Let $A(G)$ be adjacency matrix of $G$ and $\lambda_1,\lambda_2,\ldots,\lambda_n$ its eigenvalues. These are said to be the eigenvalues of the graph $G$ and to form its spectrum \cite{Cve}. The energy $E(G)$ of the graph $G$ is defined as the sum of the absolute values of its eigenvalues
$$E(G)=\sum_{i=1}^n\vert\lambda_i\vert.$$
Details and more information on graph energy can be found in \cite{Gut,Gut1,Gut2,Maj}.

\nt The Randi\'{c} matrix $R(G)=(r_{ij})_{n\times n}$ is defined as \cite{Boz,Boz1,Gut3}
\begin{displaymath}
 r_{ij}= \left\{ \begin{array}{ll}
\frac{1}{\sqrt{d_i d_j}} & \textrm{if $v_i \sim v_j$}\\
0 & \textrm{otherwise.}
\end{array} \right.
\end{displaymath}
\nt Denote the eigenvalues of the Randi\'{c} matrix $R(G)$ by $\rho_1,\rho_2,\ldots,\rho_n$ and label them in non-increasing order. The Randi\'{c} energy \cite{Boz,Boz1,Gut3} of $G$ is defined as
$$E(G)=\sum_{i=1}^n\vert\rho_i\vert.$$

\nt Two graphs $G$ and $H$ are said to be {\it Randi\'{c} energy equivalent},
or simply ${\cal RE}$-equivalent, written $G\sim H$, if
$RE(G)=RE(H)$. It is evident that the relation $\sim$ of being
${\cal RE}$-equivalence
 is an equivalence relation on the family ${\cal G}$ of graphs, and thus ${\cal G}$ is partitioned into equivalence classes,
called the {\it ${\cal RE}$-equivalence classes}. Given $G\in {\cal G}$, let
\[
[G]=\{H\in {\cal G}:H\sim G\}.
\]
We call $[G]$ the equivalence class determined by $G$.
A graph $G$ is said to be {\it Randi\'{c} energy unique}, or simply {\it ${\cal RE}$-unique}, if $[G]=\{G\}$.

\nt Similarly, we can define ${\cal E}$-equivalence for  energy and {\it ${\cal E}$-unique} for a graph.

\nt A graph $G$ is called {\it $k$-regular} if all
vertices  have the same degree $k$.  One of the famous graphs is the Petersen
graph which is a symmetric non-planar 3-regular graph.  In the study of energy and Randi\'{c} energy, it is interesting to investigate
 the characteristic polynomial and energy of this graph. We denote the Petersen graph by $P$.

 \nt In this paper, we study the energy and Randi\'{c} energy of specific graphs. In the next section, we  study energy and Randi\'c energy of  $2$-regular and $3$-regular graphs. We study cubic graphs of order 10 and list all characteristic polynomial, energy and  Randi\'{c} energy of them. As a result, we show that Petersen graph is not ${\cal RE}$-unique (${\cal E}$-unique) but can be determined by its Randi\'{c} energy (energy) and its eigenvalues. In the last section we consider some another families of graphs and study their Randi\'c characteristic polynomials.

\section{Energy of 2-regular and 3-regular graphs}

\nt The energy and Randi\'c energy of regular graphs have not been widely studied.   In this section  we consider $2$-regular and $3$-regular graphs.  The following theorem gives a relationship between the Randi\'c energy and energy of
$k$-regular graphs.

\begin{newlemma}\label{regular}\rm\cite{Rgu}
If the graph $G$ is $k$-regular then $RE(G)=\frac{1}{k}E(G)$.
\end{newlemma}

\nt Also we have the following easy lemma:

\begin{lema}\label{union}
Let $G=G_1\cup G_2\cup\ldots\cup G_m$. Then
\begin{enumerate}
\item[(i)] $E(G)=E(G_1)+ E(G_2)+\ldots+ E(G_m)$.

\item[(ii)] $RE(G)=RE(G_1)+ RE(G_2)+\ldots+ RE(G_m)$.
\end{enumerate}
\end{lema}

\nt Randi\'{c} characteristic polynomial of the cycle graph $C_n$ can be determined by the following theorem:

\begin{newlemma}\label{cycle}\rm\cite{Alikhani}
For $n\geq 3$, the Randi\'{c} characteristic polynomial of the cycle graph $C_n$ is
$$RP(C_n,\lambda)=\lambda \Lambda_{n-1}-\frac{1}{2}\Lambda_{n-2}-(\frac{1}{2})^{n-1},$$
where for every $k\geq 3$, $\Lambda_k=\lambda \Lambda_{k-1}-\frac{1}{4}\Lambda_{k-2}$ with  $\Lambda_1=\lambda$ and $\Lambda_2=\lambda ^2-\frac{1}{4}$.
\end{newlemma}

\nt By Lemma \ref{cycle}, we can find all the eigenvalues of Randi\'{c} matrix of cycle graphs. So we can compute the Randi\'{c} energy of cycles. Also every cycle is $2$-regular. By Lemma \ref{regular}, we have $E(C_n)=2RE(C_n)$. Hence we can compute energy of cycle graphs too. Every $2$-regular graph is a disjoint union of cycles. Therefore by Lemma \ref{union}, we can find energy and Randi\'{c} energy of 2-regular graphs.

%\section{Cubic graphs of order 10}

\begin{figure}[!h]
\hglue2.5cm
\includegraphics[width=11cm,height=5.1cm]{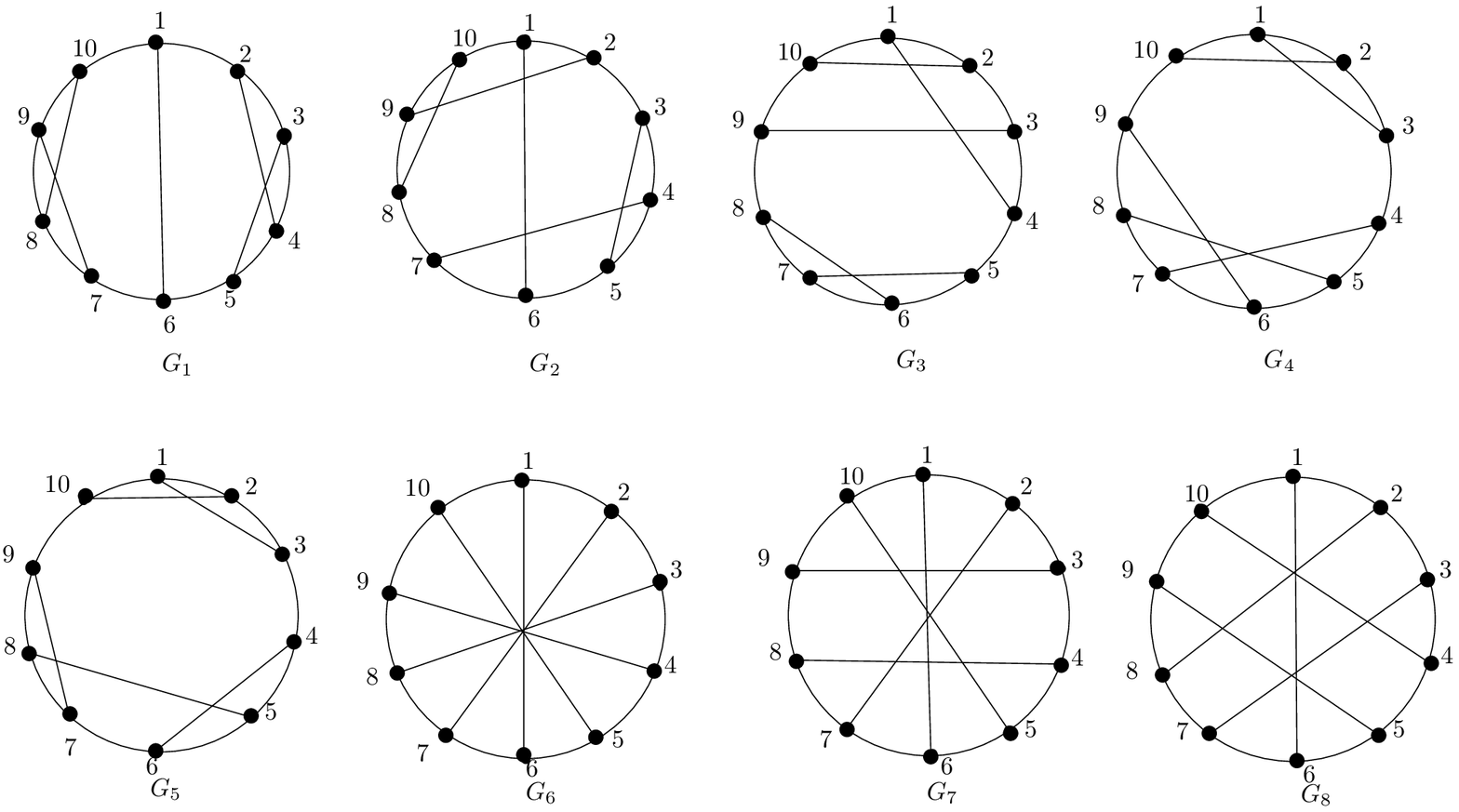}
\vglue5pt
\hglue2.5cm
\includegraphics[width=11cm,height=5.1cm]{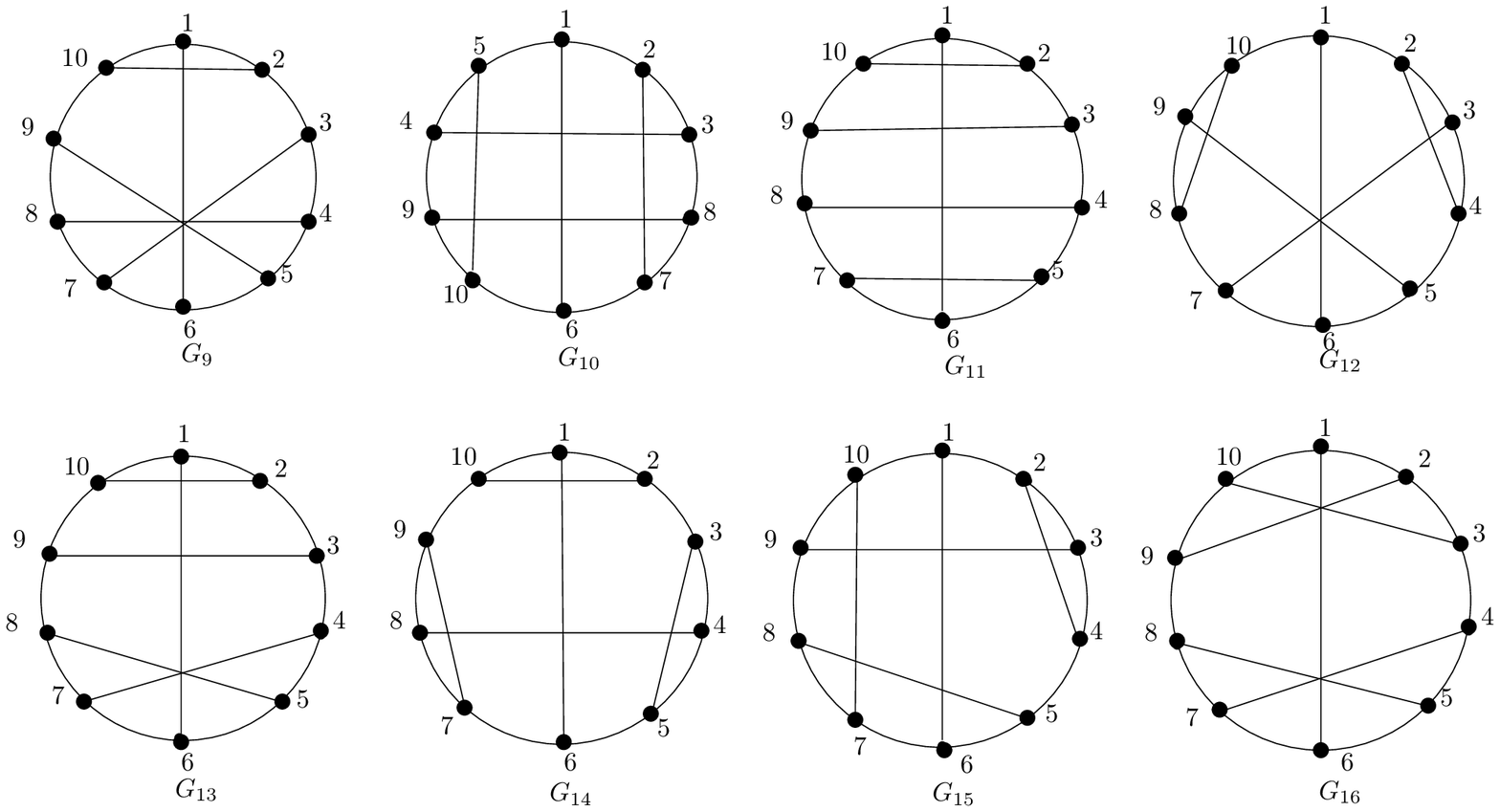}
\hglue2.5cm
\includegraphics[width=10.7cm,height=5cm]{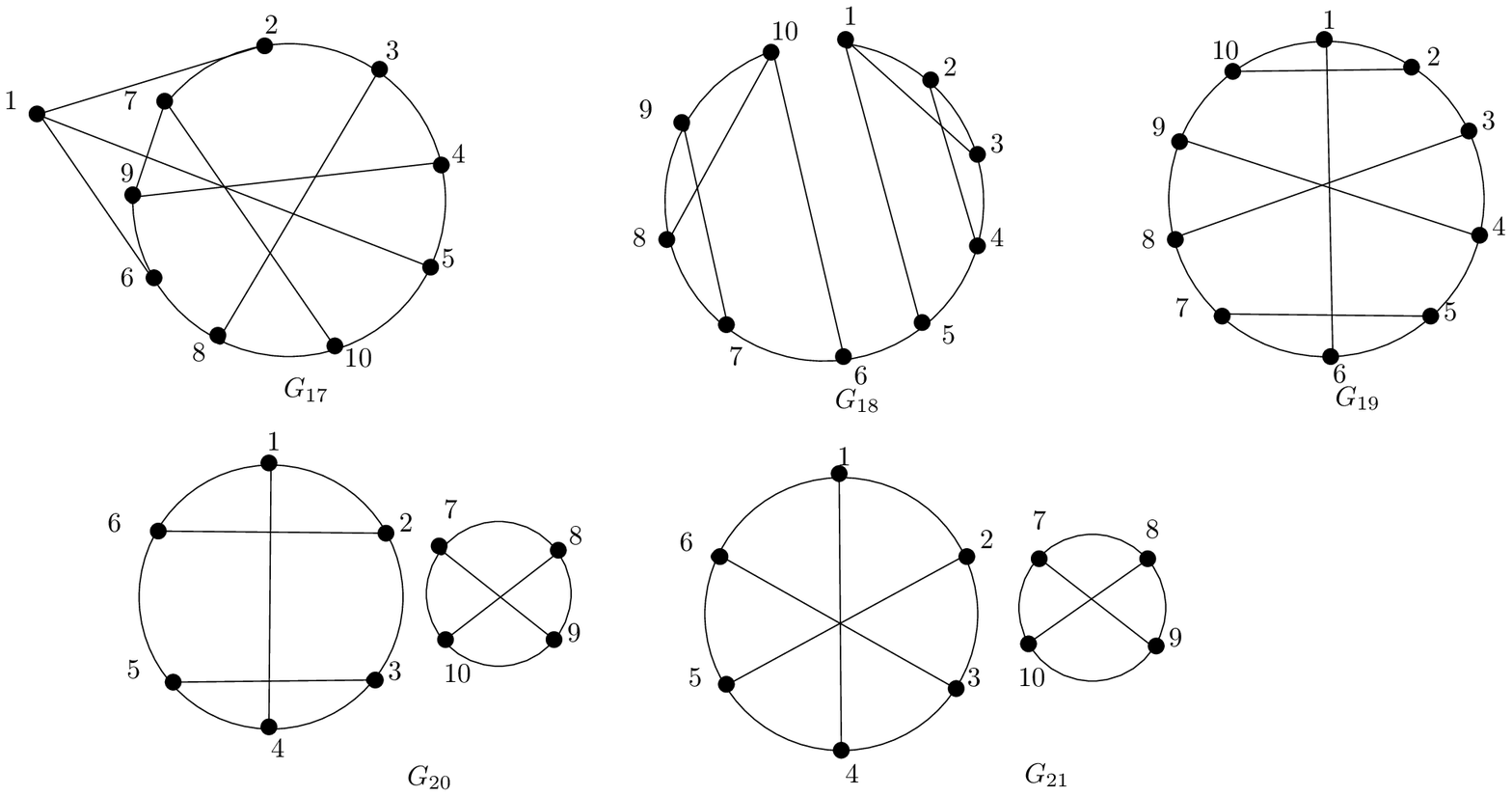}
\caption{\label{figure1} Cubic graphs of order 10. }
\end{figure}

\nt Let to consider  the characteristic polynomial of $3$-regular graphs of order $10$. Also we shall compute energy and  Randi\'{c} energy of this class of graphs. There are exactly $21$ cubic graphs of  order $10$ given in
Figure \ref{figure1} (see~\cite{reza}).

\nt We show that Petersen graph is not ${\cal RE}$-unique (${\cal E}$-unique) but can be determined by its Randi\'{c} energy (energy) and its eigenvalues. There are just two non-connected cubic graphs of order $10$. The following theorem gives us characteristic polynomial of $3$-regular graphs of order $10$. We denote the characteristic polynomial of the graph $G$ by $P(G,\lambda)$.

\nt Using Maple we computed the characteristic polynomials of $3$-regular  graphs of order $10$ in Table 1.

\bigskip
\bigskip

\begin{center}
\begin{footnotesize}
\small
\begin{tabular} {|c|c|}
\hline
$G_i$ &  $P(G_i,\lambda)$ \\
 \hline $G_1$  & $\lambda ^{10} - 15\lambda^8 - 8\lambda^7 +71\lambda^6 + 64 \lambda ^5 -101\lambda^4 -104 \lambda^3 +44\lambda^2 + 48\lambda$  \\
 \hline $G_2$ & $\lambda ^{10} - 15\lambda^8 - 4\lambda^7 +71\lambda^6 + 28 \lambda ^5 -121\lambda^4 -48 \lambda^3 +64\lambda^2 + 24\lambda$ \\
 \hline $G_3$ & $\lambda ^{10} - 15\lambda^8 - 6\lambda^7 +69\lambda^6 + 48 \lambda ^5 -96\lambda^4 -76 \lambda^3 +30\lambda^2 + 26\lambda +3$ \\
  \hline $G_4$ & $\lambda ^{10} - 14\lambda^8 - 4\lambda^7 +53\lambda^6 + 34 \lambda ^5 -48\lambda^4 -50 \lambda^3 -12\lambda^2$  \\
 \hline $G_5$  & $\lambda ^{10} - 15\lambda^8 - 8\lambda^7 +71\lambda^6 + 68 \lambda ^5 -93\lambda^4 -132 \lambda^3 -36\lambda^2$ \\
  \hline $G_6$ & $\lambda ^{10} - 15\lambda^8  +65\lambda^6  -105\lambda^4  +55\lambda^2 -9$ \\
  \hline $G_7$ &  $\lambda ^{10} - 15\lambda^8 +69\lambda^6 -12 \lambda ^5 -117\lambda^4 +36 \lambda^3 +59\lambda^2 - 12\lambda -9$\\
  \hline $G_8$ & $\lambda ^{10} - 15\lambda^8 +71\lambda^6 -16 \lambda ^5 -133\lambda^4 +64 \lambda^3 +76\lambda^2 - 48\lambda$ \\
  \hline $G_9$ & $\lambda ^{10} - 15\lambda^8 - 2\lambda^7 +71\lambda^6 + 8\lambda ^5 -132\lambda^4 -2 \lambda^3 +91\lambda^2 -8\lambda -12$ \\
  \hline $G_{10}$ & $\lambda ^{10} - 15\lambda^8  +65\lambda^6 -4\lambda ^5 -85\lambda^4 -20 \lambda^3 +35\lambda^2 +20\lambda +3$ \\
  \hline $G_{11}$ &  $\lambda ^{10} - 15\lambda^8 - 4\lambda^7 +69\lambda^6 + 32\lambda ^5 -105\lambda^4 -64 \lambda^3 +23\lambda^2 +20\lambda +3$\\
  \hline $G_{12}$ & $\lambda ^{10} - 15\lambda^8 - 4\lambda^7 +75\lambda^6 + 24\lambda ^5 -157\lambda^4 -36 \lambda^3 +144\lambda^2 +16\lambda -48$ \\
  \hline $G_{13}$ & $\lambda ^{10} - 15\lambda^8 - 2\lambda^7 +67\lambda^6 + 12\lambda ^5 -96\lambda^4 -22 \lambda^3 +35\lambda^2 +12\lambda$ \\
  \hline $G_{14}$ &  $\lambda ^{10} - 15\lambda^8  - 6\lambda^7 +75\lambda^6 + 48\lambda ^5 -144\lambda^4 -114 \lambda^3 +75\lambda^2 +68\lambda +12$\\
  \hline $G_{15}$ & $\lambda ^{10} - 15\lambda^8  - 2\lambda^7 +69\lambda^6 + 12\lambda ^5 -116\lambda^4 -24 \lambda^3 +54\lambda^2 +26\lambda +3$ \\
  \hline $G_{16}$ & $\lambda ^{10} - 15\lambda^8  +63\lambda^6  -85\lambda^4  +36\lambda^2$ \\
  \hline $G_{17}$ & $\lambda ^{10} - 15\lambda^8  +75\lambda^6 -24\lambda ^5 -165\lambda^4 +120 \lambda^3 +120\lambda^2 -160\lambda +48$ \\
  \hline $G_{18}$ & $\lambda ^{10} - 15\lambda^8  - 8\lambda^7 +63\lambda^6 + 64\lambda ^5 -37\lambda^4 -56 \lambda^3 -12\lambda^2$ \\
  \hline $G_{19}$ & $\lambda ^{10} - 15\lambda^8  - 4\lambda^7 +73\lambda^6 + 28\lambda ^5 -141\lambda^4 -52 \lambda^3 +99\lambda^2 +16\lambda -21$ \\
  \hline $G_{20}$ & $\lambda ^{10} - 15\lambda^8  - 12\lambda^7 +63\lambda^6 + 96\lambda ^5 -13\lambda^4 -84 \lambda^3 -36\lambda^2$ \\
  \hline $G_{21}$ & $\lambda ^{10} - 15\lambda^8  - 8\lambda^7 +51\lambda^6 + 72\lambda ^5 +27\lambda^4$ \\ \hline
  \end{tabular}
\end{footnotesize}
\end{center}
\begin{center}
\nt{Table 1.} Characteristic polynomial $P(G_i,\lambda)$, for $1\leq i \leq 21$.
\end{center}

\nt By computing  the roots of characteristic polynomial of cubic  graphs of order $10$, we can have the energy of these graphs.  We compute them to four decimal places. So we have table 2:

\begin{center}
\begin{footnotesize}
\small
\begin{tabular}{|c|c|c|||c|c|c|||c|c|c|} \hline
$G_i$ & $E(G_i)$ & $RE(G_i)$&$G_i$ & $E(G_i)$ & $RE(G_i)$&$G_i$ & $E(G_i)$ & $RE(G_i)$\\
\hline\hline $G_1$ & 15.1231 & 5.0410 &  $G_8$ & 15.1231 & 5.0410 & $G_{15}$ & 14.7943 &  4.9314 \\
\hline $G_2$ & 14.8596 & 4.9532 &  $G_9$ & 15.3164 & 5.1054 & $G_{16}$ & 14.0000 &  4.6666 \\
\hline $G_3$ & 14.8212 & 4.9404 &  $G_{10}$ & 14.4721 & 4.8240  & $G_{17}$ & 16.0000 & 5.3333  \\
\hline $G_4$ & 13.5143 & 4.5047 &  $G_{11}$ & 14.7020 & 4.9006 & $G_{18}$ & 13.5569 &  4.5189 \\
\hline $G_5$ & 14.2925 & 4.7641 &  $G_{12}$ & 16.0000 & 5.3333 & $G_{19}$ & 15.5791 & 5.1930  \\
\hline $G_6$ & 14.9443 & 4.9814 &  $G_{13}$ & 14.3780 & 4.7926 & $G_{20}$ & 14.0000 &  4.6666 \\
\hline $G_7$ & 15.0777 & 5.0259 &  $G_{14}$ & 15.0895 & 5.0298 & $G_{21}$ & 12.0000 & 4.0000  \\
\hline
\end{tabular}
\end{footnotesize}
\end{center}
\begin{center}
\nt{Table 2.} Energy and Randi\'c energy of cubic graphs of order $10$.
\end{center}

\begin{teorem}\label{thm3}
Six cubic graphs of order $10$   are not  ${\cal E}$-unique (${\cal RE}$-unique). If two cubic graphs of order $10$ have equal energy (Randi\'c energy), then their eigenvalues are different in exactly $3$ values.
\end{teorem}

\noindent{\bf Proof.} Using Table 2, we see that $[G_1]=\{G_1,G_8\}$, $[G_{12}]=\{G_{12},G_{17}\}$ and $[G_{16}]=\{G_{16},G_{20}\}$. Now, it suffices to find the eigenvalues of $G_1$, $G_8$, $G_{12}$, $G_{16}$ , $G_{17}$ and $G_{20}$. By Table 1 we have:
\begin{eqnarray*}
P(G_1,\lambda)&&=\lambda ^{10} - 15\lambda^8 - 8\lambda^7 +71\lambda^6 + 64 \lambda ^5 -101\lambda^4 -104 \lambda^3 +44\lambda^2 + 48\lambda \\
&&=\lambda (\lambda -3)(\lambda +2)^2 (\lambda -1)^2 (\lambda +1)^2 (\lambda -\frac{1-\sqrt{17}}{2})(\lambda -\frac{1+\sqrt{17}}{2}),
\end{eqnarray*}
\begin{eqnarray*}
P(G_8,\lambda)&&=\lambda ^{10} - 15\lambda^8 +71\lambda^6 -16 \lambda ^5 -133\lambda^4 +64 \lambda^3 +76\lambda^2 - 48\lambda\\
&&=\lambda (\lambda -3)(\lambda +2)^2 (\lambda -1)^3 (\lambda +1) (\lambda -\frac{-1+\sqrt{17}}{2})(\lambda -\frac{-1-\sqrt{17}}{2}).
\end{eqnarray*}
\nt Also
\begin{eqnarray*}
P(G_{12},\lambda)&&=\lambda ^{10} - 15\lambda^8 - 4\lambda^7 +75\lambda^6 + 24\lambda ^5 -157\lambda^4 -36 \lambda^3 +144\lambda^2 +16\lambda -48\\
&&= (\lambda -3)(\lambda -2)(\lambda +2)^3 (\lambda-1)^3 (\lambda +1)^2,
\end{eqnarray*}
\begin{eqnarray*}
P(G_{17},\lambda)&&=\lambda ^{10} - 15\lambda^8  +75\lambda^6 -24\lambda ^5 -165\lambda^4 +120 \lambda^3 +120\lambda^2 -160\lambda +48 \\
&&= (\lambda -3)(\lambda +2)^4 (\lambda-1)^5.
\end{eqnarray*}

\nt And
\begin{eqnarray*}
P(G_{16},\lambda)&&=\lambda ^{10} - 15\lambda^8  +63\lambda^6  -85\lambda^4  +36\lambda^2\\
&&= \lambda^2 (\lambda -3)(\lambda +3)(\lambda -2)(\lambda +2)(\lambda -1)^2 (\lambda +1)^2,
\end{eqnarray*}
\begin{eqnarray*}
P(G_{20},\lambda)&&= \lambda ^{10} - 15\lambda^8  - 12\lambda^7 +63\lambda^6 + 96\lambda ^5 -13\lambda^4 -84 \lambda^3 -36\lambda^2\\
&& = \lambda^2 (\lambda -3)^2(\lambda +2)^2(\lambda -1)(\lambda +1)^3.
\end{eqnarray*}
\nt So we have the result.\quad\qed

\begin{figure}[!h]
\hspace{5.5cm}
\includegraphics[width=4.5cm,height=4.5cm]{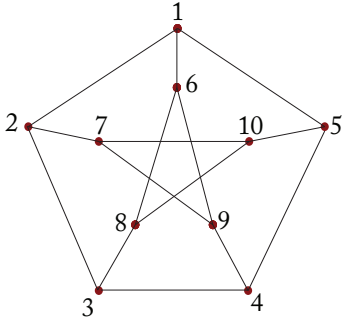}
\caption{\label{figure2} Petersen graph. }
\end{figure}

\nt Now we consider Petersen graph $P$. We have shown this graph in Figure \ref{figure2}.

\begin{teorem}
Let ${\cal G}$ be the family of $3$-regular graphs of order $10$. For the Petersen graph $P$, we have the following properties:
\begin{itemize}
\item[(i)]
$P$  is not ${\cal E}$-unique (${\cal RE}$-unique) in ${\cal G}$.
\item[(ii)]
$P$ has the maximum energy (Randi\'{c} energy) in ${\cal G}$.
\item[(iii)]
$P$ can be identify by its energy (Randi\'{c} energy) and its eigenvalues in ${\cal G}$.
\end{itemize}
\end{teorem}
\noindent{\bf Proof.}
\begin{itemize}
\item[(i)]
The adjacency matrix of $P$ is
$$A(P)=\left( \begin{array}{cccccccccc}
0&1 &0 &0 &1 &1 &0 & 0&0 & 0 \\
1&0 &1 &0 &0 &0 &1 & 0&0 & 0 \\
0&1 &0 &1 &0 &0 &0 & 1&0 & 0 \\
0&0 &1 &0 &1 &0 &0 & 0&1 & 0 \\
1&0 &0 &1 &0 &0 &0 & 0&0 & 1 \\
1&0 &0 &0 &0 &0 &0 & 1&1 & 0 \\
0&1 &0 &0 &0 &0 &0 & 0&1 & 1 \\
0&0 &1 &0 &0 &1 &0 & 0&0 & 1 \\
0&0 &0 &1 &0 &1 &1 & 0&0 & 0 \\
0&0 &0 &0 &1 &0 &1 & 1&0 & 0 \\
\end{array} \right).$$
\nt So
$det(\lambda I -A(P))=(\lambda -3)(\lambda +2)^4 (\lambda-1)^5.$
\nt Therefore we have:
$$\lambda _1=3~~,~~\lambda _2=\lambda _3=\lambda _4=\lambda _5=-2~~,~~\lambda _6=\lambda _7=\lambda _8=\lambda _9=\lambda _{10}=1,$$
\nt and so we have $E(P)=16$. By Table 2, we have $P\in \{G_{12},G_{17}\}$. Hence $P$  is not ${\cal E}$-unique (and ${\cal RE}$-unique) in ${\cal G}$.
\item[(ii)]
It follows from Part (i) and Table 2.
\item[(iii)]
It follows from Part (i) and Theorem \ref{thm3}. So $G_{17}$ is the Petersen graph.\quad\qed
\end{itemize}

\nt The following result gives a relationship between energy and permanent of adjacency matrix of two connected graphs in the family of cubic graphs of order $10$ whose have the same ${\cal E}$-equivalence class.

\begin{teorem}\label{thm5'}
If two connected cubic  graphs of order $10$ have the same energy, then their adjacency matrices have the same permanent.
\end{teorem}

\noindent{\bf Proof.} By Table 2, it suffices to find $per(A(G_1))$, $per(A(G_8))$, $per(A(G_{12}))$ and $per(A(G_{17}))$. For graph $G_1$, we have
$$A(G_1)=\left( \begin{array}{cccccccccc}
0&1 &0 &0 &0 &1 &0 & 0&0 & 1 \\
1&0 & 1& 1&0 &0 &0 &0 &0 &0  \\
0&1& 0& 1&1 &0 &0 &0 &0 &0  \\
0&1 & 1& 0&1 &0 &0 &0 &0 &0  \\
0&0 & 1& 1&0 &1 &0 &0 &0 &0  \\
1&0 & 0& 0&1 &0 &1 &0 &0 &0  \\
0&0 & 0& 0&0 &1 &0 &1 &1 &0  \\
0&0 & 0& 0&0 &0 &1 &0 &1 &1  \\
0&0 & 0& 0&0 &0 &1 &1 &0 &1  \\
1&0 & 0& 0&0 &0 &0 &1 &1 &0  \\
\end{array} \right).$$

\nt By Ryser's method, we have $per(A(G_1))=72$. Similarly we have:
$$per(A(G_8))=72 ~~,~~ per(A(G_{12}))=60 ~~,~~ per(A(G_{17}))=60.$$
So we have the result. \quad\qed

\noindent{\bf Remark 1.} The converse of Theorem \ref{thm5} is not true. Because $per(A(G_7))=per(A(G_{11}))=85$, but as we see in Table 2, $E(G_7)\neq E(G_{11})$.

\begin{korolari}
If two connected cubic  graphs of order $10$ have the same Randi\'{c} energy, then their adjacency matrices have the same permanent.
\end{korolari}

\noindent{\bf Proof.} It follows from Lemma \ref{regular}, Table 2 and Theorem \ref{thm5'}.
\quad\qed

\section{Randi\'c characteristic polynomial of a kind of Dutch-Windmill graphs}

\nt We recall that a complex number $\zeta$ is called an algebraic number (resp. an
algebraic integer) if it is a root of some monic polynomial with rational (resp.
integer) coefficients (see \cite{alg}). Since the Randi\'c characteristic  polynomial $P(G,\lambda)$ is a monic
polynomial in $\lambda$ with integer coefficients, its roots are, by definition, algebraic
integers. This naturally raises the question: Which algebraic integers can occur as
zeros of Randi\'c characteristic polynomials? And which real numbers can occur as Randi\'c energy of graphs? We are interested to numbers which are occur as Randi\'c energy. Clearly those lying in $(-\infty,2)$ are
forbidden set, because we know that if graph $G$ possesses at least
one edge, then $RE(G)\geq 2$. We think that the Randi\'c energy of graphs are dense in $[2,\infty)$. In this section  we would like to study some further results of this kind.

\begin{figure}[!h]
\hspace{1.6cm}
\includegraphics[width=11.5cm,height=4cm]{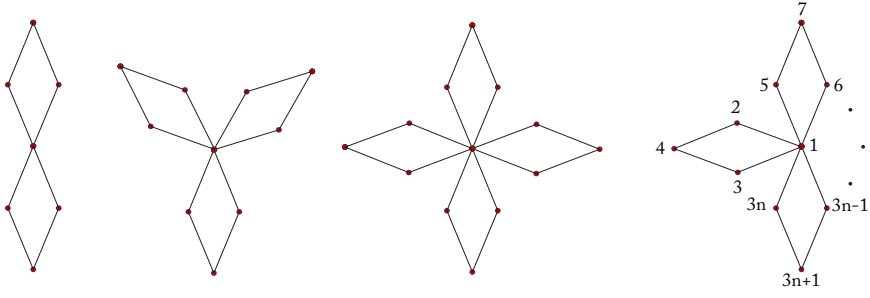}
\caption{\label{Dutch} Dutch Windmill Graph $D_4 ^2, D_4^3, D_4^4$ and $D_4^n$, respectively. }
\end{figure}

\nt Let $n$ be any positive integer and  $D_m^n$ be Dutch Windmill graph with
$(m-1)n + 1$ vertices and $mn$ edges. In other words, the  graph $D_m^n$ is a graph that can be constructed by coalescence $n$
copies of the cycle graph $C_m$ of length $m$ with a common vertex. We recall that $D_3^n$ is friendship graphs.
 Figure \ref{Dutch} shows some examples of this kind of  Dutch Windmill graphs. In this section we shall investigate the Randi\'{c} characteristic polynomial of Dutch Windmill graphs.

\nt By Lemma \ref{cycle}, we know that
$$RP(C_m,\lambda)=\lambda \Lambda_{m-1}-\frac{1}{2}\Lambda_{m-2}-(\frac{1}{2})^{m-1},$$
where for every $k\geq 3$, $\Lambda_k=\lambda \Lambda_{k-1}-\frac{1}{4}\Lambda_{k-2}$ with  $\Lambda_1=\lambda$ and $\Lambda_2=\lambda ^2-\frac{1}{4}$. We show that the Randi\'{c} characteristic polynomial of Dutch Windmill graphs can compute by the cycle which constructed it.

\begin{teorem}\label{thm5}
For $m\geq 3$, the Randi\'{c} characteristic polynomial of the  Dutch Windmill graph $D_m^n$ is
$$RP(D_m^n,\lambda)=\Lambda_{m-1} ^{n-1}.RP(C_m,\lambda),$$
where for every $k\geq 3$, $\Lambda_k=\lambda \Lambda_{k-1}-\frac{1}{4}\Lambda_{k-2}$ with  $\Lambda_1=\lambda$ and $\Lambda_2=\lambda ^2-\frac{1}{4}$.
\end{teorem}

\noindent{\bf Proof.}
\nt For every $k\geq 3$, consider

$$B_k :=
\left(\begin{array}{ccccccccc}
\lambda& \frac{-1}{2}&0&0&\ldots &0&0 &0 \\
 \frac{-1}{2}& \lambda  &\frac{-1}{2} &0&\ldots &0&0  &0 \\
0&  \frac{-1}{2}& \lambda &\frac{-1}{2} & \ldots &0&0  &0 \\
0&0&  \frac{-1}{2}& \lambda  & \ldots &0&0  &0 \\
\vdots & \vdots & \vdots & \vdots &\ddots &\vdots&\vdots &\vdots \\
 0& 0&0&0  &\ldots &\lambda&\frac{-1}{2} &0 \\
0&  0&0&0  &\ldots &\frac{-1}{2}&\lambda&\frac{-1}{2}  \\
0&  0& 0&0 &\ldots &0&\frac{-1}{2}&\lambda   \\
\end{array}\right)_{k\times k}, $$

\nt and let $\Lambda_k=det(B_k)$. It is easy to see that $\Lambda_k=\lambda \Lambda_{k-1}-\frac{1}{4}\Lambda_{k-2}$.

\nt Suppose that  $RP((D_m^n,\lambda)=det(\lambda I - R((D_m^n) )$. We have

$$ RP((D_m^n,\lambda) = det
\left(\begin{array}{c|cccc}
\lambda& A & A &\ldots &A \\
\hline
A^t &B_{m-1} &0&\ldots&0  \\
A^t  &0&B_{m-1}&\ldots&0   \\
\vdots &\vdots &\vdots&\ddots&\vdots \\
A^t  &0&0&\ldots&B_{m-1}  \\
\end{array}\right),$$

\nt where $A=\left(\begin{array}{cccccc}
\frac{-1}{2\sqrt{n}} &0&0&\ldots&0&\frac{-1}{2\sqrt{n}}
\end{array}\right)_{1\times (m-1)}$. So
$$det(\lambda I - R((D_m^n) )=\lambda\Lambda_{m-1}^n +\left(\frac{-1}{4}\Lambda_{m-2} +2( (-1)^{m+1}(\frac{-1}{2})^{m})  + (-1)^{2m+1}(\frac{1}{4})\Lambda_{m-2}\right)\Lambda_{m-1} ^{n-1}.$$

\nt Therefore

$$det(\lambda I - R((D_m^n) )=\lambda\Lambda_{m-1}^n + \left(\frac{-1}{2}\Lambda_{m-2} - (\frac{1}{2})^{m-1}\right)\Lambda_{m-1} ^{n-1}.$$

\nt Hence

$$det(\lambda I - R((D_m^n) )=\Lambda_{m-1} ^{n-1}\left(\lambda\Lambda_{m-1}-\frac{1}{2}\Lambda_{m-2} - (\frac{1}{2})^{m-1}\right)=\Lambda_{m-1} ^{n-1}RP(C_m,\lambda) . \quad\qed $$

\nt In \cite{Alikhani} we have presented two families of graphs such that their Randi\'c energy are $n+1$ and $2+(n-1)\sqrt{2}$. Here we recall the following   results:

\begin{teorem}\label{dense} {\rm \cite{Alikhani}}
\begin{enumerate}
\item[(i)]
 The Randi\'{c} energy of friendship graph $F_n$ is $RE(F_n)=n+1.$

 \item[(ii)] The Randi\'{c} energy of Dutch-Windmill graph $D_4^n$ is
$RE(D_4^n)=2+(n-1)\sqrt{2}.$
 \item[(iii)] For every $m,n\geq 2$, the Randi\'{c} energy of $K_{m,n}-e$ is $RE(K_{m,n}-e)=2+\frac{2}{\sqrt{mn}}.$
\end{enumerate}

\end{teorem}

\nt We can use Theorem \ref{thm5} to obtain $RE(D_5^n)$. Here using the definition of Randi\'c characteristic polynomial,  we prove the following result:

\begin{teorem} \label{D}
The Randi\'{c} energy of $D_5^n$ is
$$RE(D_5^n)=1+n\sqrt{5}.$$
\end{teorem}

\noindent{\bf Proof.}
The Randi\'{c} matrix of $D_5^n$ is

$$
R(D_5^n) =
\left( \begin{array}{cccccccccc}
0 & \frac{1}{2\sqrt{n}} &\frac{1}{2\sqrt{n}}&0&0&\cdots & \frac{1}{2\sqrt{n}}& \frac{1}{2\sqrt{n}}&0&0\\
\frac{1}{2\sqrt{n}}&0& 0& \frac{1}{2}&0&\ldots &0 &0&0&0 \\
\frac{1}{2\sqrt{n}}&0& 0&0& \frac{1}{2}&\ldots &0  &0&0&0\\
0&\frac{1}{2}&0 &0 &\frac{1}{2}&\ldots &0&0 &0&0\\
0&0&\frac{1}{2} &\frac{1}{2} &0&\ldots &0&0 &0&0\\
\vdots & \vdots& \vdots &\vdots & \vdots &\ddots &\vdots &\vdots  &\vdots  &\vdots\\
\frac{1}{2\sqrt{n}} & 0&0&0 &0&\ldots &0& 0& \frac{1}{2}&0  \\
\frac{1}{2\sqrt{n}} & 0&0&0 &0&\ldots &0& 0&0& \frac{1}{2}  \\
0 & 0& 0&0&0&\ldots &\frac{1}{2}&0 &0 &\frac{1}{2} \\
0 & 0& 0&0&0&\ldots &0&\frac{1}{2} &\frac{1}{2} &0\\
\end{array} \right)_{(4n+1)\times (4n+1)}.
$$

\nt Let $A= \left( \begin{array}{cccc}
 \lambda &0 &\frac{-1}{2}&0  \\
 0&\lambda &0&\frac{-1}{2} \\
\frac{-1}{2} &0 & \lambda &\frac{-1}{2}\\
0&\frac{-1}{2} &\frac{-1}{2} & \lambda \\
\end{array}\right)$
and
$C= \left( \begin{array}{cccc}
\frac{-1}{2\sqrt{n}}&0 &\frac{-1}{2}&0  \\
 \frac{-1}{2\sqrt{n}}0&\lambda &0&\frac{-1}{2} \\
0 &0 & \lambda &\frac{-1}{2}\\
0&\frac{-1}{2} &\frac{-1}{2} & \lambda \\
\end{array}\right)$.
Then

$$det(\lambda I - R((D_5^n) )=\lambda det(A)^n + \sqrt{n}det(C)det(A)^{n-1}.$$

\nt So

$$det(\lambda I - R((D_5^n) )= det(A)^{n-1} (\lambda -1)(\lambda-(\frac{\sqrt{5}}{4}-\frac{1}{4}))^2(\lambda+(\frac{\sqrt{5}}{4}+\frac{1}{4}))^2.$$

\nt Hence

$$RE(D_5^n)=1+n\sqrt{5}.\quad\qed $$

\nt Part (iii) of Theorem \ref{dense} implies that the  Randi\'c energy of graphs are dense in
$[2,3)$.  Motivated by this notation, Theorems \ref{dense} and \ref{D}, we think that  the Randi\'c energy of graphs are dense in $[2,\infty)$.   
We close this paper by the following conjecture:

\begin{konjek}
Randi\'c energy of graphs are dense in $[2,\infty)$.
\end{konjek}

\end{document}